\documentclass[11pt,reqno]{amsart}
\usepackage[breaklinks]{hyperref}
\usepackage{url}
\usepackage{amssymb}

 \renewcommand{\a}{\alpha}
\renewcommand{\b}{\beta}

\newcommand{\g}{\gamma}

\renewcommand{\(}{\left\(}
\renewcommand{\)}{\right\)}
\renewcommand{\[}{\left\[}
\renewcommand{\]}{\right\]}

\renewcommand{\i}{\infty}

\numberwithin{equation}{section}
 \theoremstyle{plain}
\newtheorem{theorem}{Theorem}[section]

\newtheorem{remark}[theorem]{Remark}
\newtheorem{corollary}[theorem]{Corollary}

\usepackage{mathrsfs}

   \makeatletter
\def\proof{\@ifnextchar[{\@oproof}{\@nproof}}
\def\@oproof[#1][#2]{\trivlist\item[\hskip\labelsep\textit{#2 Proof of\
#1.}~]\ignorespaces}
\def\@nproof{\trivlist\item[\hskip\labelsep\textit{Proof.}~]\ignorespaces}

\makeatother

\newmuskip\pFqmuskip

\newcommand*\pFq[6][8]{%
  \begingroup 
  \pFqmuskip=#1mu\relax
  \mathchardef\normalcomma=\mathcode`,
  \mathcode`\,=\string"8000
  \begingroup\lccode`\~=`\,
  \lowercase{\endgroup\let~}\pFqcomma
  {}_{#2}F_{#3}{\left[\genfrac..{0pt}{}{#4}{#5};#6\right]}%
  \endgroup
}
\newcommand{\pFqcomma}{{\normalcomma}\mskip\pFqmuskip}

\begin{document}
\title[Extended Higher Herglotz function \textup{II} ]{Extended Higher Herglotz function \textup{II}}

\author{Rajat Gupta and Rahul Kumar}
\address{Discipline of Mathematics, Indian Institute of Technology, Gandhinagar, Palaj, Gandhinagar 382355, Gujarat, India}\email{rajat\_gupta@iitgn.ac.in}

\address{Center for Geometry and Physics, Institute for Basic Science (IBS), Pohang 37673, Republic of Korea}\email{rahul@ibs.re.kr}

\maketitle

\begin{abstract}
Very recently, Radchenko and Zagier revived the theory of Herglotz functions. The main goal of the article is to show that one of the formulas on page 220 of Ramanujan's Lost Notebook actually lives in the realms of this theory. As a consequence of our general theorem, we derive an interesting identity analogous to Ramanujan's formula for $\zeta(2m+1)$. We also introduce a character analogue of the Herglotz function and initiate its theory by obtaining an elegant functional equation governed by it.
\end{abstract}

\section{Introduction}

In his seminal work on the Kronecker limit formula for a real quadratic field, Zagier \cite{zagier} found an interesting function which is now known as the Herglotz function. It is defined by
\begin{align}\label{rd herglotz defn}
F(x):=\sum_{n=1}^\infty\frac{\psi(nx)-\log(nx)}{n},\qquad x\in\mathbb{C}\backslash(-\infty,0],
\end{align}
where $\psi(x):=\Gamma'(x)/\Gamma(x)$ is the digamma function. Prior to Zagier, a function quite similar to  $F(x)$ is also appeared in the work of Herglotz \cite{herglotz}. This is why Radchenko and Zagier call it the \emph{Herglotz function} in \cite{rz}. Zagier establishes beautiful functional equations for $F(x)$, namely, \cite[Equations (7.4), (7.8)]{zagier}, for $x\in\mathbb{C}\backslash(-\infty,0]$,
{\allowdisplaybreaks\begin{align}
F(x)-F(x+1)-F\left(\frac{x}{x+1}\right)&=-F(1)+\textup{Li}_2\left(\frac{1}{1+x}\right),\nonumber\\
F(x)+F\left(\frac{1}{x}\right)&=2F(1)+\frac{1}{2}\log^{2}(x)-\frac{\pi^2}{6x}(x-1)^2,\label{fe2}
\end{align}}
where $\mathrm{Li}_2(x):=\sum_{n=1}^\infty t^n/n^2$ is dilogarithm function and \cite[Equation (7.12)]{zagier}
\begin{equation*}
F(1)=-\frac{1}{2}\g^2-\frac{\pi^2}{12}-\g_1.
\end{equation*}
Here $\gamma$ and $\g_1$ denote the Euler and first Stieltjes constants respectively. In the same paper \cite[Section 8]{zagier}, Zagier used these functional equations to prove Meyer's theorem \cite{meyer}. 

Recently, Radchenko and Zagier \cite{rz}  studied $F(x)$ extensively and unearthed a myriad of connections that $F(x)$ has with other areas of number theory. They provided many exciting properties of $F(x)$, including functional equations, special values at rational or quadratic irrational arguments,  connections with Stark's conjecture and with 1-cocycles for the modular group PSL$(2,\mathbb{Z})$.



Many generalizations of $F(x)$ or its variants have been studied in the literature. Vlasenko and Zagier \cite{vz},  in their work on higher Kronecker limit formula, considered the study of the higher Herglotz function:
\begin{align}\label{vz herglotz}
F_k(x):=\sum_{n=1}^\infty\frac{\psi(nx)}{n^k}, \qquad k\in\mathbb{N},\ k>1,\ x\in\mathbb{C}\backslash(-\infty,0].
\end{align}
They found the following functional equations for $F_k(x)$ \cite[Equations (11), (12)]{vz}:
{\allowdisplaybreaks\begin{align}\label{vz1f}
F_k(x)+(-x)^{k-1}F_k\left(\frac{1}{x}\right)&=-\gamma\zeta(k)\left(1+(-x)^{k-1}\right)-\sum_{r=2}^{k-1}\zeta(r)\zeta(k+1-r)(-x)^{r-1}\nonumber\\
&\quad+\zeta(k+1)\left((-x)^{k}-\frac{1}{x}\right),
\end{align}}
and\footnote{The interpretation of $\zeta(1,k-1)$ given in \cite[p.~28]{vz} has a minus sign missing in front of the whole term, that is, the correct form is $-\left(\zeta(k-1,1)+\zeta(k)-\g\zeta(k-1)\right)$.} 
\begin{align}
&F_k(x)-F_k(x+1)+(-x)^{k-1}F_k\left(\frac{x+1}{x}\right)\nonumber\\
&=(-x)^{k-1}\left(\zeta(k,1)+\zeta(k+1)-\g\zeta(k)\right)-\sum_{r=1}^{k-1}\zeta(k+1-r,r)(-x)^{r-1}\nonumber\\
&\hspace{3.5cm}+\zeta(k+1)\left(\frac{(-x)^{k}}{x+1}-\frac{1}{x}\right),\nonumber
\end{align}
where, $\zeta(k)$ is the Riemann zeta function and
\begin{equation*}
\zeta(m, n):=\sum_{p>q>0}\frac{1}{p^mq^n}\hspace{5mm}(m\geq2, n\geq1).
\end{equation*}
Vlasenko and Zagier \cite[p.~54, Section 3]{vz} also dealt with another twisted extension of \eqref{vz herglotz} to treat the zeta functions of ray classes. Ishibashi studied $j^{th}$ order Herglotz function to find an explicit representation of the Laurent series coefficients of the zeta function associated to indefinite quadratic forms \cite[Theorem 3]{ishibashi}.

Let $\alpha:\mathbb{Z}\to\mathbb{C}$ be a periodic function with period $M$. Masri \cite{masri} defined the following $L$-series generalization of $F(x)$: 
\begin{align*}
F(\alpha,s,x):=\sum_{n=1}^\infty\frac{\alpha(n)\left(\psi(nx)-\log(nx)\right)}{n^s},\qquad\mathrm{Re}(s)>0,\ x>0,
\end{align*} 
and provided its analytic continuation and evaluated it at integer points.

Very recently, the authors along with Dixit \cite{dgk} studied the \emph{extended higher Herglotz function}
\begin{align*}
\mathscr{F}_{k,N}(x):=\sum_{n=1}^{\infty}\frac{\psi(n^{N}x)-\log(n^{N}x)}{n^k},\qquad x\in\mathbb{C}\backslash(-\infty,0],
\end{align*}
where $k$ and $N$ are positive real numbers such that $k+N>1$. Note that the function $\mathscr{F}_{k,N}(x)$ subsumes Herglotz function $F(x)$ as well as Vlasenko and Zagier's higher Herglotz function $F_k(x)$. One of the nice functional equations which $\mathscr{F}_{k,N}(x)$ satisfies is {\allowdisplaybreaks\begin{align}\label{hhfeeqn}
&x^{\frac{1-k}{N}}\mathscr{F}_{k,N}(x)-\frac{(-1)^k}{N}\sum_{j=-(N-1)}^{(N-1)}{\vphantom{\sum}}''e^{i\pi j (k-1)/N}\mathscr{F}_{\frac{N+k-1}{N},\frac{1}{N}}\left(\frac{e^{-\frac{i\pi j}{N}}}{x^{1/N}}\right)\nonumber\\
&=\frac{\pi}{N} \frac{\zeta\left(1+\frac{k-1}{N}\right)}{\sin \left(\frac{\pi}{N}(k-1) \right)}+x^{\frac{1-k}{N}}\left(-\left(\gamma+\log x \right) \zeta(k)+N\zeta'(k)\right)-\frac{1}{x^{\frac{N+k-1}{N}}}\zeta(k+N)+\mathscr{B}(k, N, x),
\end{align}}
where  $k,N\in \mathbb{N}$ such that $1<k\leq N$ and $x\in\mathbb{C}\backslash(-\infty,0]$ and
\begin{equation*}
\mathscr{B}(k, N, x):=\left\{
	\begin{array}{ll}
		(-1)^{k+N+1}x^{1/N}\zeta\left( 1+\frac{k}{N}\right) & \mbox{if } k=N \\
		0 & \mbox{if } k\neq N
	\end{array}
\right. .
\end{equation*}
Here and throughout this paper, the notation $\sum_{j=-(N-1)}^{(N-1)}{\vphantom{\sum}}''$ denotes a sum over over $j=-(N-1), -(N-3),\cdots, N-3, N-1$. For many other nice properties of $\mathscr{F}_{k,N}(x)$ and non-trivial applications of its functional equations, we refer the reader to \cite{dgk}. 

In the last section of their paper, the present authors with Dixit \cite{dgk} posed several questions. Here in this article, we answer some of those questions affirmatively. Before we do this though, we record an elegant formula from Ramanujan's Lost Notebook on page 220.
If  
\begin{align*}
\varphi(x):=\psi(x)+\frac{1}{2x}-\log x,
\end{align*}
then for $\a,\b$ positive such that $\a\b=1$, we have \cite[p.~220]{lnb}
\begin{align}\label{ramanujan formula}
\sqrt{\a}&\left\{\frac{\gamma - \log(2 \pi \a)}{2\a} +\sum_{n=1}^{\infty}\varphi(n\a)\right\}=\sqrt{\b}\left\{\frac{\gamma - \log(2 \pi \b)}{2\b} +\sum_{n=1}^{\infty}\varphi(n\b)\right\}.
\end{align}
He also provided the following surprising integral representation which is equal to either sides of \eqref{ramanujan formula}:
\begin{align}\label{riemann c integral}
-\frac{1}{\pi^{3/2}}\int_{0}^{\infty}\left|\Xi\left( \frac{t}{2}\right) \Gamma\left(\frac{-1+it}{4} \right)\right|^2\frac{\cos\left(\frac{t}{2}\log\a \right)}{1+t^2}dt,
\end{align}
where $\Xi(t)$ denotes Riemann's function defined by \cite[p.~16, Equations (2.1.14), (2.1.12)]{titch}
\begin{align*}
\Xi(t)&:=\xi(\tfrac{1}{2}+it),
\end{align*}
and
\begin{align*}
\xi(s)&:=\frac{s}{2}(s-1)\pi^{-\frac{1}{2}s}\Gamma(\tfrac{1}{2}s)\zeta(s).
\end{align*}
The proof of the above formula can be found in \cite{bd}. The above transformation and various other transformations of its type have generated interesting mathematics. We refer the reader to recent papers \cite{dg}, \cite{dk}, and references therein for the extensive literature related to this formula.

It was commented in \cite[p.~27]{dgk} that ``\emph{it seems that \eqref{ramanujan formula} lives in the realms of the theory of Herglotz function}" by observing the similarity between the definition  of $F(x)$ in \eqref{rd herglotz defn} and series appearing on both sides of \eqref{ramanujan formula}. One of the main objectives of this paper is to show that this is in-fact the case, that is, Ramanujan's formula \eqref{ramanujan formula} is actually a special case of two of our more general results on the Herglotz functions (Theorems \ref{N=1 and general k}, and \ref{general N and k=0}).

Note that the functional equation \eqref{hhfeeqn} is valid only when $1<k\leq N$. It was also asked in \cite[Section 5, p.~28]{dgk} if there exists a functional equation for $\mathscr{F}_{k,N}(x)$ where $k$ is greater than $N$? We answer this question affirmatively in Theorem \ref{general theorem}.

Although $\mathscr{F}_{k,N}(x)$ is a generalization of Valsenko-Zagier's higher Herglotz function $\mathscr{F}_{k}(x)$, the functional equation of $\mathscr{F}_{k}(x)$ in \eqref{vz1f} could not be obtained from \eqref{hhfeeqn} as a special case due to the restriction $1<k\leq N$. This problem is now circumvented through our Theorem \ref{general theorem}. 

To achieve our goals, we work with the following version of the extended higher Herglotz function:
\begin{align}\label{hhhf}
\mathfrak{F}_{k,N}(x):=\sum_{n=1}^\infty\frac{1}{n^k}\left(\psi\left(n^Nx\right)-\log\left(n^Nx\right)+\frac{1}{2n^Nx}\right),
\end{align}
for $x\in\mathbb{C}\backslash(-\infty,0],\ \mathrm{and}\ k+2N>1$. The series is absolutely and uniformly convergent for $k+2N>1$ since \cite[p.~259, formula 6.3.18]{abr} $$\psi(x)\sim\log(x)-\frac{1}{2x}-\frac{1}{12x^2}+\cdots,\qquad\mathrm{as}\ x\to\infty,\ |\mathrm{arg}(x)|<\pi.$$

Note that for $k=N=1$, $\mathfrak{F}_{k,N}(x)$ reduces to the function $F^*(x)$:
\begin{align*}
F^*(x):=\sum_{n=1}^\infty\frac{1}{n}\left(\psi\left(nx\right)-\log\left(nx\right)+\frac{1}{2nx}\right).
\end{align*}
The function $F^*(x)$ was used by Radchenko and Zagier \cite[p.~17, Section 7.3]{rz} to reveal the cocycle nature of the Herglotz function $F(x)$ by finding an Eichler-type integral for it.

The Herglotz function $F(x)$ is a special case of $\mathfrak{F}_{k,N}(x)$:
\begin{align*}
\mathfrak{F}_{1,1}(x)=F(x)+\frac{\pi^2}{12x}.
\end{align*}
Also,
\begin{align*}
\mathfrak{F}_{k,N}(x)=\mathscr{F}_{k,N}(x)+\frac{1}{2x}\zeta(k+N).
\end{align*}

\begin{remark}
Ramanujan's formula \eqref{ramanujan formula} can be rephrased in terms of $\mathfrak{F}_{k,1}(x)$ by
\begin{align}\label{k1}
\mathfrak{F}_{k,1}(x)-\frac{1}{x}\mathfrak{F}_{k,1}\left(\frac{1}{x}\right)=\frac{1}{2}\left(\gamma-\log\left(\frac{2\pi}{x}\right)\right)-\frac{1}{2x}\left(\gamma-\log(2\pi x)\right).
\end{align}
Thus Ramanujan understood the importance of such a function. Therefore, it would not be unfair to say that Ramanujan studied Hergltoz-type functions even before Herglotz and Zagier. Further, he went beyond the functional equation to find a surprising intergal representation for this Herglotz-type function which involves the Riemann $\Xi$-function.

\end{remark}



The next objective of this paper is to initiate the theory of a character analogue of the Herglotz function. We define it by
\begin{align}\label{Fxi}
F_{k}(x,\chi):=\sum_{n=1}^{\i}\frac{\chi(n)\psi(nx,\chi)}{n^k},\quad x\in\mathbb{C}\backslash(-\infty,0]
\end{align}
here and throughout the article $\chi(n)$ is a primitive, nonprincipal character modulo $d$, $k$ is a non-negative integer, and
\begin{align*}
\psi(x,\chi):=-\sum_{k=1}^{\i}\frac{\chi(k)}{k+x}.
\end{align*}
For a real character $\chi$, the function $\psi(n,\chi)$ is equivalent to the character analogue
of the digamma function obtained by the logarithmic differentiation of the Weierstrass product form of the character analogue of the gamma function for real
characters introduced by Berndt \cite{b}:
\begin{align*}
\Gamma(s,\chi):=e^{-sL(1,\chi)}\prod_{n=1}^{\i}\left(1+\frac{s}{n}\right)^{-\chi(n)}e^{s\chi(n)/n},\hspace{5mm} s\in \mathbb{C}\backslash \mathbb{Z}_{<0}.
\end{align*}

Note that the series in \eqref{Fxi} is convergent for any $k\geq0$ since \cite[p.~334, Corollary 4.4]{ad}
\begin{align*}
\psi(x,\chi)\sim-\frac{L(0,\chi)}{x}-\chi(-1)\sum_{j=2}^\infty \frac{B_j(\overline{\chi})}{jx^j},\quad |\mathrm{arg}(x)|<\pi,
\end{align*}
as $x\to\infty$. Here $B_j(\chi)$ is the generalized Bernoulli numbers \cite[p.~426]{b}.

\section{Main Results}\
We divide this section into two parts. The first part is devoted to extending the modular relation of Ramanujan \eqref{ramanujan formula}, and the other is dedicated to initiating the theory of $F_k(x,\chi)$ defined in \eqref{Fxi}.

\subsection{Extension of the formula on page 220 of Ramanujan's Lost Notebook}

In this subsection, we state three of our main theorems. The first is Theorem \ref{N=1 and general k}, which is valid for $k\in\mathbb{N}\cup\{0\}$ and $N=1$. The second is Theorem \ref{general N and k=0} which holds for $k=1$ and $N\in\mathbb{N}$. We conclude this subsection by providing Theorem \ref{general theorem} which holds for any natural numbers $k$ and $N$ greater than 1. 
\begin{theorem}\label{N=1 and general k}
Let $\mathfrak{F}_{k,N}(x)$ be defined in \eqref{hhhf}. For any  $k\in\mathbb{N}\cup\{0\}$ and $\mathrm{Re} (x)>0$, the following identity holds:
\begin{align*}
\mathfrak{F}_{k,1}(x)=(-1)^kx^{k-1}\left\{\mathfrak{F}_{k,1}\left(\frac{1}{x}\right)-\mathcal{B}_k(x)\right\},
\end{align*}
where, for $k=0$,
\begin{align*}
\mathcal{B}_k(x):=\frac{x}{2}\left(\log\left(\frac{2\pi}{x}\right)-\gamma\right)+\frac{1}{2}\left(\gamma-\log(2\pi x)\right),
\end{align*}
for $k=1$,
\begin{align*}
\mathcal{B}_k(x):=\frac{1}{6}\left(3\log^2(x)+\pi^2-6\gamma^2-12\gamma_1\right)-\frac{\pi^2}{12}\left(x+\frac{1}{x}\right),
\end{align*}
and for $k\geq2$,
\begin{align*}
\mathcal{B}_k(x):&=\zeta(k)\left(\log(x)-\gamma\right)+\zeta'(k)+(-1)^kx^{1-k}\left((\log(x)+\gamma)\zeta(k)-\zeta'(k)\right)+\frac{(-1)^k}{2x^k}\zeta(k+1)\nonumber\\
&\qquad-\frac{x}{2}\zeta(1+k)+\sum_{j=1}^{k-2}(-1)^{j-1}\zeta(j+1)\zeta(k-j)x^{-j}.
\end{align*}
\end{theorem}

\begin{remark}
It is easy to see that Ramanujan's formula \eqref{ramanujan formula}\textup{(}or \eqref{k1}\textup{)}, Zagier's functional equation \eqref{fe2}, and Vlasenko-Zagier's functional equation \eqref{vz1f} are special cases of the above theorem for $k=0,~k=-1,$ and $k\geq2$, respectively.
\end{remark}

An immediate consequence of the above theorem is the following beautiful symmetric formula.
\begin{theorem}\label{alphabeta form}
Let $k$ be any natural number greater than 1. Let $\mathrm{Re}(\alpha)>0$ and $\mathrm{Re}(\beta)>0$ such that $\alpha\beta=1$, then the following relation holds:
\begin{align}\label{alphabeta form eqn}
&\alpha^{\frac{1-k}{2}}\left\{(\gamma+\log(\alpha))\zeta(k)-\zeta'(k)+\frac{1}{2\alpha}\zeta(k+1)+\mathfrak{F}_{k,1}(\alpha)\right\}\nonumber\\
&=(-1)^k\beta^{\frac{1-k}{2}}\left\{(\gamma+\log(\beta))\zeta(k)-\zeta'(k)+\frac{1}{2\beta}\zeta(k+1)+\mathfrak{F}_{k,1}(\beta)\right\}\nonumber\\
&\qquad+\sum_{j=1}^{k-2}(-1)^j\zeta(j+1)\zeta(k-j)\alpha^{\frac{1-k+j}{2}}\beta^{-\frac{j}{2}}.
\end{align}
\end{theorem}
This formula is analogous to the Ramanujan's formula for $\zeta(2m+1)$ \cite[p.~173, Ch. 14, Entry 21(i)]{ramnote}, namely, for $\textup{Re}(\a)>0, \textup{Re}(\b)>0, \a\b=\pi^2$ and $m\in\mathbb{Z}\backslash\{0\}$, by
\begin{align*}
&\a^{-m}\left\{\frac{1}{2}\zeta(2m+1)+\sum_{n=1}^{\infty}\frac{n^{-2m-1}}{e^{2\a n}-1}\right\}\nonumber\\
&=(-\b)^{-m}\left\{\frac{1}{2}\zeta(2m+1)+\sum_{n=1}^{\infty}\frac{n^{-2m-1}}{e^{2\b n}-1}\right\}-2^{2m}\sum_{j=0}^{m+1}\frac{(-1)^jB_{2j}B_{2m+2-2j}}{(2j)!(2m+2-2j)!}\a^{m+1-j}\b^j,
\end{align*}
where $B_n$ is the $n^{\textup{th}}$ Bernoulli number.

The next theorem is another extension of \eqref{ramanujan formula} in the second parameter $N$.

\begin{theorem}\label{general N and k=0}
Let $N$ be any natural number and $\mathrm{Re} (x)>0$. Then
\begin{align}\label{general N and k=0 eqn}
\mathfrak{F}_{0,N}(x)=\frac{x^{-1/N}}{N}\sideset{}{''}\sum_{j=-(N-1)}^{(N-1)}e^{-\frac{i\pi j}{N}}\mathfrak{F}_{\frac{N-1}{N},\frac{1}{N}}\left(\frac{e^{-\frac{i\pi j}{N}}}{x^{1/N}}\right)-\frac{1}{2}\left(\log\left(\frac{(2\pi)^N}{x}\right)-\gamma\right)+R(N,x),
\end{align}
where
\begin{align}\label{rnx}
R(N,x):=
\begin{cases}
\frac{1}{2x}\left(\log(2\pi x)-\gamma\right), &\mathrm{if} \ N=1\\
-\frac{\pi x^{-1/N}\zeta\left(1-1/N\right)}{N\sin(\pi/N)}-\frac{1}{2x}\zeta(N), &\mathrm{if} \ N>1.
\end{cases}
\end{align}
\end{theorem}
Now we state a general theorem which is valid for any $k, N\in \mathbb{N}\backslash{\{1\}}.$
\begin{theorem}\label{general theorem}
Let $k,\ N$ be any natural numbers greater $1$. Then for $\mathrm{Re} (x)>0$, we have
\begin{align*}
\mathfrak{F}_{k,N}(x)=\frac{(-1)^k}{N}x^{\frac{k-1}{N}}\sideset{}{''}\sum_{j=-(N-1)}^{N-1}e^{\frac{i\pi j(k-1)}{N}}\mathfrak{F}_{\frac{N+k-1}{N},\frac{1}{N}}\left(\frac{e^{-\frac{i\pi j}{N}}}{x^{1/N}}\right)-\frac{x^{\frac{k-1}{N}}}{N}\mathscr{R}_{N,k}(x),
\end{align*}
where, for $N\nmid (k-1)$ and $N\nmid k$,
\begin{align*}
\mathscr{R}_{N,k}(x)&=\frac{N}{2}x^{\frac{1-k-N}{N}}\zeta(k+N)-\frac{\pi\zeta\left(\frac{N+k-1}{N}\right)}{\sin\left(\frac{\pi}{N}(k-1)\right)}-Nx^{\frac{1-k}{N}}\left(N\zeta'(k)-(\gamma+\log x)\zeta(k)\right)\nonumber\\
&\qquad +N\sum_{j=1}^{\lfloor\frac{k-1}{N}\rfloor}(-1)^j\zeta(1+j)\zeta(k-N)x^{\frac{k-1}{N}+j},
\end{align*}
and, for $N|k$, 
\begin{align*}
\mathscr{R}_{N,k}(x)&=\frac{N}{2}x^{\frac{1-k-N}{N}}\zeta(k+N)-\frac{\pi\zeta\left(\frac{N+k-1}{N}\right)}{\sin\left(\frac{\pi}{N}(k-1)\right)}-Nx^{\frac{1-k}{N}}\left(N\zeta'(k)-(\gamma+\log x)\zeta(k)\right)\nonumber\\
&\quad-\frac{N}{2}(-1)^{\frac{k}{N}}x^{\frac{1}{N}}\zeta\left(\frac{k+N}{N}\right) +N\sum_{j=1}^{\lfloor\frac{k-1}{N}\rfloor}(-1)^j\zeta(1+j)\zeta(k-N)x^{\frac{k-1}{N}+j},
\end{align*}
and, for $N|(k-1)$,
\begin{align*}
\mathscr{R}_{N,k}(x)&=\frac{N}{2}x^{\frac{1-k-N}{N}}\zeta(k+N)+(-1)^{\frac{k-1}{N}}\left((N\gamma-\log x)\zeta\left(\frac{N+k-1}{N}\right)-\zeta'\left(\frac{N+k-1}{N}\right)\right)\nonumber\\
&\quad -Nx^{\frac{1-k}{N}}\left(N\zeta'(k)-(\gamma+\log x)\zeta(k)\right)+N\sum_{j=1}^{\lfloor\frac{k-1}{N}\rfloor-1}(-1)^j\zeta(1+j)\zeta(k-N)x^{\frac{k-1}{N}+j}.
\end{align*}
\end{theorem}

When we restrict $1<k\leq N$ in the above theorem,  it gives \eqref{hhfeeqn} as a special case:
\begin{corollary}
The functional equation \eqref{hhfeeqn} holds true.
\end{corollary}

\subsection{Character analogue of the Herglotz function}
For a non-principal Dirichlet character $\chi(n)$, let $L(s, \chi)$ denote the Dirichlet $L$-function defined by
\begin{align*}
L(s,\chi)=\sum_{n=1}^{\i}\frac{\chi(n)}{n^s}
\end{align*}
for Re$(s)>1$. This series converges conditionally for $0<\textup{Re}(s)<1$. Also, it can be analytically continued to an entire function of $s$. 

The two-term functional equation governed by $F_k(x,\chi)$ is given in the following theorem.
\begin{theorem}\label{2.7}
Let $F_{k}(x,\chi)$ be defined in \eqref{Fxi}. Let $\chi(n)$ be a primitive, nonprincipal character modulo $d$. Then $k \in \mathbb{N}\cup\{0\}$ and \textup{Re}$(x)>0$, we have
\begin{align*}
F_{k}(x,\chi)-(-1)^kx^{k-1}F_{k}\left(x,\frac{1}{\chi}\right)=\sum_{j=0}^{k-1}(-x)^{k-j-1}L(1+j,\chi)L(-j+k, \chi).
\end{align*}
\end{theorem}

Now to state our next theorem, we define
\begin{align*}
b:=\left\{
	\begin{array}{ll}
		0  & \mbox{if } \chi(-1)=1 \\
		1 & \mbox{if } \chi(-1)=-1.
	\end{array}
\right.
\end{align*}
Then the $\Xi(s,\chi)$ can be defined as
\begin{align}\label{Xi}
\Xi(t,\chi):=\xi\left(\frac{1}{2}+it,\chi\right),
\end{align}
where
\begin{align}\label{xi}
\xi(s,\chi):=\left( \frac{\pi}{d}\right)^{-(s+b)/2}\Gamma\left(\frac{s+b}{2} \right)L(s,\chi).
\end{align}

Our next theorem gives an elegant modular relation along with the integral comprising $\Xi(t,\chi)$-function.
\begin{theorem}\label{2.8}
Let $\chi(n)$ be a primitive, nonprincipal character modulo $d$. Further assume $\mathrm{Re}(\alpha)>0,\ \mathrm{Re}(\beta)>0$ such that $\a\b=1$, then
\begin{align*}
\sqrt{\a}F_{0}(\a,\chi)&=\sqrt{\b}F_{0}(\b,\chi)\nonumber\\
&=-\frac{1}{2\pi}\left(\frac{\pi}{d}\right)^{b+\frac{1}{2}}\int_{0}^{\i}\left|\frac{\Gamma\left(\frac{1}{2}+\frac{it}{2} \right)}{\Gamma\left(\frac{1}{4}+\frac{it}{4}+\frac{b}{2} \right)}\right|^2\Xi\left(-\frac{t}{2},\chi\right)\Xi\left(\frac{t}{2},\chi\right)\cos\left( \frac{t}{2}\log\a\right)dt.
\end{align*}
\end{theorem}

The above theorem reduces to  \cite[Corollaries 5.1 and 5.2]{ad} depending on the parity of the character.


\section{Proofs}

We first present the proof of our general Theorem \ref{general theorem}. Throughout the proofs, the notation $\int_{(c)}ds$  will denote the line integral $\int_{c-i\infty}^{c+i\infty}ds$ with $c=\mathrm{Re}(s)$.
\begin{proof}[Theorem \text{\ref{general theorem}}][]
Kloosterman's formula is given by \cite[p.~25, Equation (2.9.2)]{titch}
\begin{align}\label{tich}
\psi(x+1)-\log(x)=\frac{1}{2\pi i}\int_{(c)}\frac{-\pi\zeta(1-z)}{\sin(\pi z)}x^{-z}dz,
\end{align}
which is valid for $0<c=\mathrm{Re}(z)<1$. An application of the functional equation 
\begin{align*}
\psi(x+1)=\psi(x)+1/x
\end{align*}
in \eqref{tich} implies that
\begin{align}\label{abc}
\psi(x)-\log(x)+\frac{1}{x}=\frac{1}{2\pi i}\int_{(c)}\frac{-\pi\zeta(1-z)}{\sin(\pi z)}x^{-z}dz.
\end{align}
For our purpose we need the integral representation for the function $\psi(x)-\log(x)+\frac{1}{2x}$. For that, we shift the line of the integration to $1<d=\mathrm{Re}(z)<2$. Let $\mathfrak{C}$ be a positively oriented rectangular contour formed by the points $c-iT, d-iT, d+iT, c+iT$. Observe that the integrand has only one simple pole at $s=1$ due to $\sin(\pi z)$ and its residue is $-1/(2x)$ as $\zeta(0)=-1/2$. Therefore, residue theorem immediately yields
\begin{align}\label{T}
\frac{1}{2\pi i}\left(\int_{d-iT}^{d+iT}+\int_{d+iT}^{c+iT}-\int_{c-iT}^{c+iT}+\int_{c-iT}^{d-iT}\right)\frac{-\pi\zeta(1-z)}{\sin(\pi z)}x^{-z}dz=-\frac{1}{2x}.
\end{align}
Note that the integrals along the horizontal lines goes to zero by using the Stirling formula in the vertical strip $p\leq\sigma\leq q$ \cite[p.~224]{cop}:
\begin{equation}\label{strivert}
  |\Gamma(s)|=\sqrt{2\pi}|t|^{\sigma-\frac{1}{2}}e^{-\frac{1}{2}\pi |t|}\left(1+O\left(\frac{1}{|t|}\right)\right)
\end{equation}
as $|t|\to \infty$. Thus, from \eqref{T}, we have 
\begin{align}\label{T2}
\frac{1}{2\pi i}\int_{(d)}\frac{-\pi\zeta(1-z)}{\sin(\pi z)}x^{-z}dz=\frac{1}{2\pi i}\int_{(c)}\frac{-\pi\zeta(1-z)}{\sin(\pi z)}x^{-z}dz-\frac{1}{2x}.
\end{align}
Combining \eqref{abc} and \eqref{T2} together, we obtain 
\begin{align}\label{klooster}
\psi(x)-\log(x)+\frac{1}{2x}=\frac{1}{2\pi i}\int_{(d)}\frac{-\pi\zeta(1-z)}{\sin(\pi z)}x^{-z}dz.
\end{align}
Replacing $x$ by $n^Nx$ in the above equation and invoking \eqref{hhhf}, we obtain
\begin{align}\label{before change of variable}
\mathfrak{F}_{k,N}(x)&=\sum_{n=1}^\infty\frac{1}{n^k}\frac{1}{2\pi i}\int_{(d)}\frac{-\pi\zeta(1-z)}{\sin(\pi z)}(n^Nx)^{-z}dz\nonumber\\
&=\frac{1}{2\pi i}\int_{(d)}\frac{-\pi\zeta(1-z)\zeta(k+Nz)}{\sin(\pi z)}x^{-z}dz,
\end{align}
where in the last step we interchanged the order of the summation and integration and used the series definition of $\zeta(s)$. We now make change of variable $s=1-k-Nz$ in \eqref{before change of variable} to deduce, for $1-k-2N<c<1-k-N$,
\begin{align}\label{before cauchy}
\mathfrak{F}_{k,N}(x)&=\frac{x^{\frac{k-1}{N}}}{N}\frac{1}{2\pi i}\int_{(c)}\frac{-\pi\zeta\left(1+\frac{s+k-1}{N}\right)\zeta(1-s)}{\sin\left(\frac{\pi}{N}(1-s-k)\right)}x^{\frac{s}{N}}ds
\end{align}
We next shift the line of integration to $1<\lambda<2$. Consider the contour formed by the line segments $[\lambda-iT,\lambda+iT],\ [\lambda+iT,c+iT],\ [c+iT,c-iT]$ and $[c-iT,\lambda-iT]$. In order to do so, we encounter several poles of different orders depending upon $k$ and $N$. We denote the sum of the residues by $\mathscr{R}_{k,N}(x)$. By invoking Stirling's formula \eqref{strivert}, we see that the integrals along the horizontal lines vanish as $T\to\infty$. Thus, by residue theorem and \eqref{before cauchy}, we have
\begin{align}\label{cauchy}
\mathfrak{F}_{k,N}(x)&=\frac{x^{\frac{k-1}{N}}}{N}\left(\frac{1}{2\pi i}\int_{(\lambda)}\frac{\pi\zeta\left(1+\frac{s+k-1}{N}\right)\zeta(1-s)}{\sin\left(\frac{\pi}{N}(s+k-1)\right)}x^{\frac{s}{N}}ds-\mathscr{R}_{k,N}(x)\right).
\end{align}
Now, the main task is to evaluate the line integral present on the right-hand side of the above equation. Let us denote it by $\mathcal{I}_{k,N}(x)$. To that end, we invoke following result \cite[Lemma 4.1]{dgkm}
\begin{align*}
\frac{1}{\sin(z)}=\frac{1}{\sin(Nz)}\sideset{}{''}\sum_{j=-(N-1)}^{N-1}e^{ijz}
\end{align*}
with $z=\frac{\pi}{N}(s+k-1)$ in $\mathcal{I}_{k,N}(x)$, so that
\begin{align*}
\mathcal{I}_{k,N}(x)=\sum_{j=-(N-1)}^{(N-1)}{\vphantom{\sum}}''e^{\frac{i\pi j (k-1)}{N}}\frac{1}{2\pi i}\int_{(\lambda)}\frac{\pi\zeta(1-s)\zeta\left(1+\frac{s+k-1}{N}\right)}{\sin \left(\pi(s+k-1) \right)}\left(\frac{e^{-\frac{i\pi j}{N}}}{x^{1/N}}\right)^{-s}.
\end{align*}
We can use the series definition of $\zeta\left(1+\frac{s+k-1}{N}\right)$ since $1<\lambda=\mathrm{Re}(s)<2$. Thus, we arrive at
\begin{align*}
\mathcal{I}_{k,N}(x)=(-1)^{k}\sum_{j=-(N-1)}^{(N-1)}{\vphantom{\sum}}''e^{\frac{i\pi j (k-1)}{N}}\sum_{n=1}^{\i}\frac{1}{n^{1+\frac{k-1}{N}}}\frac{1}{2\pi i}\int_{(\lambda)}\frac{-\pi\zeta(1-s)}{\sin \left(\pi s\right)} \left(\left(\frac{n}{x}\right)^{\frac{1}{N}}e^{-\frac{i\pi j}{N}}\right)^{-s}ds.
\end{align*}
Invoking \eqref{klooster} in the above equation so as to obtain
\begin{align}\label{line integral evaluation}
\mathcal{I}_{k,N}(x)&=(-1)^{k}\sum_{j=-(N-1)}^{(N-1)}{\vphantom{\sum}}''e^{\frac{i\pi j (k-1)}{N}}\sum_{n=1}^{\i}\frac{1}{n^{1+\frac{k-1}{N}}}\left\{\psi\left(\left(\frac{n}{x}\right)^{\frac{1}{N}}e^{-\frac{i\pi j}{N}}\right)-\log\left(\left(\frac{n}{x}\right)^{\frac{1}{N}}e^{-\frac{i\pi j}{N}}\right)\right.\nonumber\\
&\left.\qquad\qquad\qquad+\frac{1}{2}\left(\frac{n}{x}\right)^{-\frac{1}{N}}e^{\frac{i\pi j}{N}}\right\}\nonumber\\
&=(-1)^{k}\sum_{j=-(N-1)}^{(N-1)}{\vphantom{\sum}}''e^{\frac{i\pi j (k-1)}{N}}\mathfrak{F}_{\frac{N+k-1}{N},\frac{1}{N}}\left(\frac{e^{-\frac{i\pi j}{N}}}{x^{1/N}}\right),
\end{align}
where in the last step we employed the definition of $\mathfrak{F}_{k,N}(x)$ from \eqref{hhhf}.

We next calculate the residues. We break this calculation into different cases depending on whether $N$ divides $k$ and $k-1$ or not. It boils down to only three cases: $N\nmid k$ and $N\nmid(k-1)$, $N|k$ and $N|(k-1)$ as if $N$ divides $k$ then it cannot divide $k-1$ and vice-versa. \\

\textbf{Case I:} Let $N\nmid k$ and $N\nmid (k-1)$. Observe that the integrand has simple poles at $s=0$ due to $\zeta(1-s),$ and at $s=1-k-N, 1-k+Nj$, where $1\leq j\leq \lfloor\frac{k-1}{N}\rfloor$, due to  $\sin\left(\frac{\pi}{N}(1-s-k)\right)$. It also has double pole at $s=1-k$ due to $\sin\left(\frac{\pi}{N}(1-s-k)\right)$ and $\zeta\left(1+\frac{s+k-1}{N}\right)$. These can be evaluated to as:\\
\begin{align}\label{not devinding both}
R_{0}&=-\frac{\pi\zeta\left(\frac{N+k-1}{N}\right)}{\sin\left(\frac{\pi}{N}(k-1)\right)}\nonumber\\
R_{1-k-N}&=\frac{N}{2}\zeta(k+N)x^{\frac{1-k-N}{N}}\nonumber\\
R_{1-k}&=N\left\{N\zeta'(k)-(\gamma+\log(x))\zeta(k)\right\}x^{\frac{1-k}{N}}\nonumber\\
R_{1-k+Nj}&=N\sum_{j=1}^{\lfloor\frac{k-1}{N}\rfloor}(-1)^j\zeta(1+j)\zeta(k-Nj)x^{\frac{k-1}{N}+j}.
\end{align}
\textbf{Case II:} Let $N|k$. Therefore, $k=mN$ for some $m\in\mathbb{N}$. This implies $\sin\left(\frac{\pi}{N}(1-s-k)\right)=(-1)^{m+1}\sin\left(\frac{\pi}{N}(1-s)\right)$.  Hence, the integrand has an extra simple pole at $s=1$ besides all other poles in the previous case. The residue at the pole $s=1$ is:
\begin{align}\label{deviding k}
R_1=-\frac{N}{2}(-1)^{k/N}x^{1/N}\zeta\left(\frac{k+N}{N}\right).
\end{align}
\textbf{Case III:} Let $N|(k-1)$. In this case, $\sin\left(\frac{\pi}{N}(1-s-k)\right)$ also has pole at $s=0$ which was not the case before (this was $j=\lfloor\frac{k-1}{N}\rfloor$ case in Case I). This implies that the integrand has now double pole at $s=0$ along with other poles in the  Case I, but now $j$ runs between $1\leq j\leq \lfloor\frac{k-1}{N}\rfloor-1$. Therefore, we only need to calculate the residue at $s=0$, that is:
\begin{align}\label{deviding k-1}
R_0=(-1)^{\frac{k-1}{N}}\left\{\big(N\gamma-\log(x)\big)\zeta\left(\frac{N+k-1}{N}\right)-\zeta'\left(\frac{N+k-1}{N}\right)\right\}.
\end{align}
We now combine all of the residues from \eqref{not devinding both}, \eqref{deviding k}, and \eqref{deviding k-1} together, then it is nothing but the term $\mathscr{R}_{k,N}(x)$ defined in the statement of the theorem according to the different cases.

Theorem now follows upon using the above facts and \eqref{line integral evaluation} in \eqref{cauchy}.
\end{proof}

\begin{proof}[Theorem \text{\ref{N=1 and general k}}][]
Note that up to \eqref{before cauchy}, the calculation holds for any $N\in\mathbb{N}$ and $k\in\mathbb{N}\cup\{0\}$. Thus, letting $N=1$ in \eqref{before cauchy} gives for $-1-k<c<-k$,
\begin{align*}
\mathfrak{F}_{k,1}(x)=(-1)^kx^{k-1}\frac{1}{2\pi i}\int_{(c)}\frac{-\pi\zeta(s+k)\zeta(1-s)}{\sin(\pi s)}x^sds.
\end{align*}
We shift the line of integration to $1<\lambda<2$ in the similar manner as we did in after \eqref{before cauchy}. Thus, we have
\begin{align}\label{cauchy 2}
\mathfrak{F}_{k,1}(x)=(-1)^kx^{k-1}\left(\frac{1}{2\pi i}\int_{(\lambda)}\frac{-\pi\zeta(s+k)\zeta(1-s)}{\sin(\pi s)}x^sds-\mathscr{R}_{k,1}(x)\right).
\end{align}
The calculation of the residues depends on the values of $k$ is $0$, $1$ and $k\geq2$. \\
\textbf{Case I:} When $k=0$, the integrand has double pole at $s=0$ due to $\zeta(1-s)$ and $\sin(\pi s)$ and $s=1$ due to $\zeta(s)$ and $\sin(\pi s)$. The residues at these points are:
\begin{align}\label{k=0}
R_0&=\frac{1}{2}\big(\gamma-\log(2\pi x)\big)\nonumber\\
R_1&=\frac{x}{2}\left(\log\left(\frac{2\pi}{x}\right)-\gamma\right).
\end{align}
\textbf{Case II:} Let $k=1$. Now integrand has simple pole at $s=1,-1$ due to $\sin(\pi s)$ and pole of order three at $s=0$ due to $\sin(\pi s), \zeta(s+1)$ and $\zeta(1-s)$ with residues:
\begin{align}\label{k=1}
R_1&=-\frac{\pi^2}{12}x\nonumber\\
R_{-1}&=\frac{\pi^2}{12}\frac{1}{x}\nonumber\\
R_0&=\frac{1}{6}\left(3\log^2(x)+\pi^2-6\gamma^2-12\gamma_1\right).
\end{align}
\textbf{Case III:} For $k\geq2$, the integrand has pole of order 2 at $s=0$ and $s=1-k$ because of $\sin(\pi s), \zeta(1-s)$ and $\sin(\pi s),\zeta(s+k)$, respectively. It has also simple poles at $s=1,-k$ and $s=-j, 1\leq j\leq k-2$ due to $\sin(\pi s)$. The residues at these poles can be evaluated to 
\begin{align}\label{k>2}
R_0&=\zeta(k)(\log(x)-\gamma)+\zeta'(k)\nonumber\\
R_{1-k}&=(-1)^kx^{1-k}\left\{\zeta(k)(\log(x)+\gamma)-\zeta'(k)\right\}\nonumber\\
R_1&=-\frac{x}{2}\zeta(1+k)\nonumber\\
R_{-k}&=\frac{(-1)^k}{2x^k}\zeta(k+1)\nonumber\\
R_j&=\sum_{j=1}^{k-2}(-1)^{j-1}\zeta(j+1)\zeta(k-j)x^{-j}.
\end{align}
Upon combining the residues from \eqref{k=0}, \eqref{k=1} and \eqref{k>2} together in \eqref{cauchy 2}, we see that $\mathscr{R}_{k,1}(x)$ is nothing but $\mathcal{B}_k(x)$ defined in the statement of the theorem.

The same argument, we used to obtain \eqref{line integral evaluation}, can be adapted to show that
\begin{align*}
\frac{1}{2\pi i}\int_{(\lambda)}\frac{-\pi\zeta(s+k)\zeta(1-s)}{\sin(\pi s)}x^sds=\mathfrak{F}_{k,1}\left(\frac{1}{x}\right).
\end{align*}
Substituting the residues and the above integral in \eqref{cauchy 2}, we complete the proof.
\end{proof}

\begin{proof}[Theorem \text{\ref{general N and k=0}}][]
The argument is similar to the proof of Theorem \ref{general theorem}. One can proceed along the same lines with letting $k=0$ and then only difference is while calculating the residual terms. When $N>1$, the integrand has simple poles at $s=0$, $s=1-N$ and a pole of order two at $s=1$. Whereas, for $N=1$, it has poles of order two at $s=0$ and $s=1$. Therefore, we have
\begin{align}\label{c}
\mathfrak{F}_{0,N}(x)&=\frac{x^{-\frac{1}{N}}}{N}\left(\frac{1}{2\pi i}\int_{(\lambda)}\frac{\pi\zeta\left(1+\frac{s-1}{N}\right)\zeta(1-s)}{\sin\left(\frac{\pi}{N}(s-1)\right)}x^{\frac{s}{N}}ds-\mathscr{R}_{0,N}(x)\right).
\end{align}
It can be easily seen that the residual term $\mathscr{R}_{0,N}$ present in above equation is
\begin{align*}
-\frac{N}{2}x^{1/N}\left(\log\left(\frac{(2\pi)^N}{x}\right)-\gamma\right)+R(N,x),
\end{align*}
where $R(N,x)$ is defined in \eqref{rnx}. 

The line integral in \eqref{c} can be evaluated by letting $k=0$ in \eqref{line integral evaluation}. Hence, now \eqref{general N and k=0 eqn} follows straightforward upon using the above facts.
\end{proof}

\begin{proof}[Theorem \text{\ref{alphabeta form}}][]
For $k\geq2$, Theorem \ref{N=1 and general k} gives
\begin{align}
x^{\frac{1-k}{2}}\mathfrak{F}_{k,1}(x)&=(-1)^kx^{\frac{k-1}{2}}\mathfrak{F}_{k,1}\left(\frac{1}{x}\right)-x^{\frac{1-k}{2}}\left((\log x+\gamma)\zeta(k)-\zeta'(k)\right)\nonumber\\
&\quad+(-1)^kx^{\frac{k-1}{2}}\left(\left(\log\left(\frac{1}{x}\right)+\gamma\right)\zeta(k)-\zeta'(k)\right)-\frac{1}{2}x^{-\frac{1+k}{2}}\zeta(k+1)\nonumber\\
&\quad+\frac{1}{2}(-1)^kx^{\frac{k+1}{2}}\zeta(k+1)+\sum_{j=1}^{k-2}(-1)^{k+j-1}\zeta(j+1)\zeta(k-j)x^{\frac{k-1}{2}-j}.\nonumber
\end{align}
Let $x=\alpha$ and $\alpha\beta=1$. Then the above equation can be rephrased as
\begin{align}\label{zx}
&\alpha^{\frac{1-k}{2}}\left\{(\log\alpha+\gamma)\zeta(k)-\zeta'(k)+\frac{1}{2\alpha}\zeta(k+1)+\mathfrak{F}_{k,1}(\alpha)\right\}\nonumber\\
&=(-1)^k\beta^{\frac{1-k}{2}}\left\{(\log\beta+\gamma)\zeta(k)-\zeta'(k)+\frac{1}{2\beta}\zeta(k+1)+\mathfrak{F}_{k,1}(\beta)\right\}\nonumber\\
&\qquad+\sum_{j=1}^{k-2}(-1)^{k+j-1}\zeta(j+1)\zeta(k-j)\alpha^{\frac{k-1}{2}-j}.
\end{align}
Replacing $j$ by $k-1-j$ in the finite sum of \eqref{zx} and using the fact $\alpha\beta=1$, we arrive at \eqref{alphabeta form eqn}.
\end{proof}

We now present the proofs of our results on the character analogue of the Herglotz function $F_k(x,\chi)$. We first prove the two-term functional equation.
\begin{proof}[Theorem \text{\ref{2.7}}][]
We first evaluate the sum
\begin{align*}
\sum_{n=1}^{\i}\frac{\chi(n)\psi(nx,\chi)}{n^k},
\end{align*}
by employing \cite[Corollary 4.2]{ad}\footnote{It is easy to see that this result is valid for Re$(x)>0.$}, hence we have for $0<c<1$,  
\begin{align*}
\sum_{n=1}^{\i}\frac{\chi(n)\psi(nx,\chi)}{n^k}&=\sum_{n=1}^{\i}\frac{\chi(n)}{n^k}\frac{1}{2\pi i}\int_{(c)}\frac{-\pi L(1-s,\chi)}{\sin(\pi s)}(nx)^{-s}\ ds\\
&=\frac{1}{2\pi i}\int_{(c)}\frac{-\pi L(1-s,\chi)L(k+s,\chi)}{\sin(\pi s)}x^{-s}\ ds.
\end{align*}
Upon performing the change of variable from $s$ to $1-s-k$, we obtain
\begin{align}\label{xisi}
\sum_{n=1}^{\i}\frac{\chi(n)\psi(nx,\chi)}{n^k}
&=\frac{(-1)^k}{2\pi i}\int_{(d)}\frac{-\pi L(s+k,\chi)L(1-s,\chi)}{\sin(\pi s)}x^{s+k-1}\ ds,
\end{align}
where $-k<d<1-k$. To evaluate it further, we shift the line of the integration to $0<c=\mathrm{Re}(s)<1$. Let $\mathfrak{C}$ be a positively oriented rectangular contour formed by the points $c-iT, d-iT, d+iT, c+iT$. Observe that the integrand has simple poles at $s=0, 1,\cdots, 1-k$ due to $\sin(\pi s)$. Therefore, residue theorem immediately yields
\begin{align*}
\frac{1}{2\pi i}\left(\int_{c-iT}^{c+iT}+\int_{c+iT}^{d+iT}-\int_{d-iT}^{d+iT}+\int_{d-iT}^{c-iT}\right)\frac{(-1)^{k+1} x^{k-1}\pi L(1-s,\chi)L(s+k,\chi)}{\sin(\pi s)}x^{s}\ ds=-\sum_{j=0}^{k-1}R_{-j}.
\end{align*}
Note that, as $T \to \i$ the contribution from integrals along the vertical lines is zero. Hence, 
\begin{align*}
\frac{(-1)^{k+1} x^{k-1}}{2\pi i}&\int_{(d)}\frac{\pi L(1-s,\chi)L(1+s,\chi)}{\sin(\pi s)}x^{s}\ ds\nonumber\\
&=\frac{(-1)^{k+1} x^{k-1}}{2\pi i}\int_{(c)}\frac{\pi L(1-s,\chi)L(1+s,\chi)}{\sin(\pi s)}x^{s}\ ds-\sum_{j=0}^{k-1}R_{-j}, 
\end{align*}
where,
\begin{align*}
R_{-j}=(-1)^{1+j-k}x^{-1-j+k}L(1+j,\chi)L(-j+k, \chi).
\end{align*}
Finally, from \eqref{xisi} and invoking the definition \eqref{Fxi} we complete our proof.
\end{proof}

\begin{proof}[Theorem \text{\ref{2.8}}][]
The first-equality follows by letting $k=0$ in Theorem \ref{2.7}. Now, for the second equality we evaluate
\begin{align*}
\sum_{n=1}^{\i}\chi(n)\psi(nx,\chi).
\end{align*}
Hence, for $0<c<1,$ from \eqref{xisi} and by employing the reflection formula for the gamma function $\Gamma(s)\Gamma(1-s)=\pi/\sin(\pi s)$ in the second step, we see that
\begin{align*}
\sum_{n=1}^{\i}\chi(n)\psi(nx,\chi)&=\sum_{n=1}^{\i}\chi(n)\frac{1}{2\pi i}\int_{(c)}\frac{-\pi L(1-s,\chi)}{\sin(\pi s)}(nx)^{-s}\ ds\\
&=-\frac{1}{2\pi i}\int_{(c)}\Gamma(s)\Gamma(1-s) L(1-s,\chi)L(s,\chi)x^{-s}\ ds.
\end{align*}
Applying the definition \eqref{xi} in the above equation to obtain
\begin{align*}
\sum_{n=1}^{\i}\chi(n)\psi(nx,\chi)
&=-\frac{1}{2\pi i}\int_{(c)}\frac{\Gamma(s)\Gamma(1-s)}{\Gamma\left(\frac{1-s+b}{2} \right)}\left(\frac{\pi}{d} \right)^{\frac{1-s+b}{2}}\frac{\xi(1-s,\chi)\xi
(s,\chi)}{\Gamma\left(\frac{s+b}{2} \right)}\left(\frac{\pi}{q} \right)^{\frac{s+b}{2}}x^{-s}\ ds.
\end{align*}
Letting $c=\frac{1}{2},$ we see that 
\begin{align*}
\sum_{n=1}^{\i}\chi(n)\psi(nx,\chi)
&=-\left(\frac{\pi}{d} \right)^{\frac{1}{2}+b}\frac{1}{2\pi i}\int_{\frac{1}{2}-i\i}^{\frac{1}{2}+i\i}\frac{\Gamma(s)\Gamma(1-s)}{\Gamma\left(\frac{1-s+b}{2} \right)\Gamma\left(\frac{s+b}{2} \right)}\xi(1-s,\chi)\xi(s,\chi)x^{-s}\ ds.
\end{align*}
Now we perform the change of variable $s=\frac{1}{2}+it$ and employing \eqref{Xi}, 
\begin{align*}
\sum_{n=1}^{\i}\chi(n)\psi(nx,\chi)
&=-\frac{1}{\sqrt{x}}\left(\frac{\pi}{d} \right)^{\frac{1}{2}+b}\frac{1}{2\pi}\int_{-\i}^{\i}\frac{\Gamma\left(\frac{1}{2}+it \right)\Gamma\left(\frac{1}{2}-it \right)}{\Gamma\left(\frac{1}{4}-\frac{it}{2}+\frac{b}{2} \right)\Gamma\left(\frac{1}{4}+\frac{it}{2}+\frac{b}{2} \right)}\Xi(-t,\chi)\Xi(t,\chi)x^{-it}\ dt.
\end{align*}
We perform the change of variable by $t\to t/2$, and after some considerable simplification, we arrive at 
{\allowdisplaybreaks\begin{align*}
&\sum_{n=1}^{\i}\chi(n)\psi(nx,\chi)\nonumber\\
&=-\frac{1}{2\pi\sqrt{x}}\left(\frac{\pi}{d} \right)^{b+\frac{1}{2}}\int_{0}^{\i}\frac{\Gamma\left(\frac{1}{2}+\frac{it}{2}  \right)\Gamma\left(\frac{1}{2}-\frac{it}{2} \right)}{\Gamma\left(\frac{1}{4}-\frac{it}{4}+\frac{b}{2} \right)\Gamma\left(\frac{1}{4}+\frac{it}{4}+\frac{b}{2} \right)}\Xi\left(-\frac{t}{2} ,\chi\right)\Xi\left(\frac{t}{2} ,\chi\right)\cos\left(\frac{t \log x}{2}\right)dt.
\end{align*}}
This completes the proof of the theorem.
\end{proof}

\section{Acknowledgements}
The authors would like to show their sincere gratitude to Prof. Atul Dixit for fruitful discussions and suggestions on the manuscript. The first author’s research was supported by the SERB-DST CRG grant CRG/2020/002367 of Prof. Atul Dixit and partly by his institute IIT Gandhinagar. The second author’s research was supported by the grant IBS-R003-D1 of the IBS-CGP, POSTECH, South Korea. Both authors sincerely thank their respective funders for the support.

\end{document}